\documentclass[11pt]{amsart}
\usepackage{amsmath,amssymb,amsthm}
\usepackage[latin1]{inputenc}
\usepackage{version,tabularx,multicol}
\usepackage{graphicx,float}

\newcommand{\reff}[1]{(\ref{#1})}

\theoremstyle{plain}
\newtheorem{theo}{Theorem}[section]

\newtheorem{prop}[theo]{Proposition}

\newtheorem{lem}[theo]{Lemma}

\theoremstyle{remark}
\newtheorem{rem}[theo]{Remark}

\newcommand{\cc}{{\mathcal C}}

\newcommand{\cf}{{\mathcal F}}

\newcommand{\cj}{{\mathcal J}}

\newcommand{\cn}{{\mathcal N}}
\newcommand{\cm}{{\mathcal M}}

\newcommand{\cs}{{\mathcal S}}

\newcommand{\E}{{\mathbb E}}

\newcommand{\N}{{\mathbb N}}
\renewcommand{\P}{{\mathbb P}}

\newcommand{\R}{{\mathbb R}}

\newcommand{\ind}{{\bf 1}}

\newcommand{\supp}{{\rm supp}\;}

\newcommand{\Card}{{\rm Card}\;}

\newcommand{\inv}[1]{\mathop{\frac{1}{ #1}}\nolimits}
\newcommand{\expp}[1]{\mathop {\mathrm{e}^{ #1}}}

\begin{document}

\title[Asymptotics for the small fragments]{Asymptotics for the small fragments of the fragmentation at nodes}

\date{\today}
\author{Romain Abraham} 

\address{
MAPMO, F\'ed\'eration Denis Poisson, Université d'Orléans,
B.P. 6759,
45067 Orléans cedex 2
FRANCE
}
  
\email{romain.abraham@univ-orleans.fr}

\author{Jean-François Delmas}

\address{ENPC-CERMICS, 6-8 av. Blaise Pascal,
  Champs-sur-Marne, 77455 Marne La Vallée, France.}

\email{delmas@cermics.enpc.fr}

\begin{abstract}
  We consider the  fragmentation at nodes of the L\'evy continuous random tree
  introduced  in a  previous paper.   In this  framework we  compute the
  asymptotic for  the number of  small fragments at time  $\theta$. This
  limit   is  increasing   in   $\theta$  and   discontinuous.  In   the
  $\alpha$-stable  case  the fragmentation  is  self-similar with  index
  $1/\alpha$,  with $\alpha  \in (1,2)$ and the  results are  close to
  those Bertoin obtained  for general self-similar fragmentations but
  with an additional assumtion which is not fulfilled here.
\end{abstract}



\maketitle

\section{Introduction}

A  fragmentation process  is a  Markov  process which  describes how  an
object with  given total  mass  evolves  as it breaks  into several
fragments randomly as time passes.  Notice there may be loss of mass but
no  creation.  Those processes  have been  widely studied  in the
recent  years, see Bertoin~\cite{b:rfcp} and  references therein.  To be
more precise, the  state space of a fragmentation process  is the set of
the non-increasing sequences of masses with finite total mass:
\[
\mathcal{S}^{\downarrow}=\left\{s=(s_1,s_2,\ldots); \;    s_1\ge    s_2\ge
  \cdots\ge                                         0\quad\text{and}\quad
  \Sigma(s)=\sum_{k=1}^{+\infty}s_k<+\infty\right\}.
\]
If we  denote by  $P_s$ the law  of a  $\cs^{\downarrow}$-valued process
$\Lambda       =(\Lambda       (\theta),\theta\ge       0)$      starting       at
$s=(s_1,s_2,\ldots)\in\cs^{\downarrow}$,  we say  that $\Lambda  $  is a
fragmentation  process if  it is  a Markov  process such  that $\theta\mapsto
\Sigma(\Lambda(\theta))$   is   non-increasing  and   if   it  fulfills   the
fragmentation property: the law  of $(\Lambda(\theta),\theta\ge 0)$ under $P_s$ is
the non-increasing reordering of  the fragments of independent processes
of respective  laws $P_{(s_1,0,\ldots)}$,$P_{(s_2,0,\ldots)}$, \ldots. In
other words,  each fragment  after dislocation behaves  independently of
the others,  and its evolution depends  only on its initial  mass.  As a
consequence, to describe  the law of the fragmentation  process with any
initial condition, it suffices to study the laws $P_r:=P_{(r,0,\ldots)}$
for any  $r\in (0,+\infty)$, i.e.  the law of the  fragmentation process
starting with a single mass $r$.

A fragmentation process is said to be self-similar of index $\alpha'$ if,
for  any $r>0$,  the process  $\Lambda$
under $P_r$  is distributed as the  process $(r\Lambda (r^{\alpha'} \theta),\theta\ge
0)$  under $P_1$.  Bertoin~\cite{b:ssf}  proved that  the law  of a
self-similar   fragmentation   is  characterized   by:   the  index   of
self-similarity $\alpha'$,  an erosion coefficient $c$ which  corresponds to a
 rate of  mass loss,  and a  dislocation measure  $\nu$ on
$\cs^{\downarrow}$ which describes sudden  dislocations of a fragment of
mass~1. The dislocation measure of a fragment of size $r$, $\nu_r$ is given by
$\int F(s) \nu_r(ds)= r^{\alpha'} \int F(rs) \nu(ds)$.

When  there  is   no  loss  of  mass  (which   implies  that  $c=0$  and
$\alpha'>0$),  under   some  additional  assumptions,   the  number  of
fragments at a  fixed time is infinite. A  natural question is therefore
to study  the asymptotic behavior when  $\varepsilon$ goes down  to 0 of
$N^{\varepsilon}(\theta)=\Card\{i,       \Lambda_i(\theta)>\varepsilon\}$
where  $\Lambda(\theta)=(\Lambda_1(\theta),\Lambda_2(\theta),\ldots)$ is
the    state   of    the   fragmentation    at   time    $\theta$,   see
Bertoin~\cite{b:smssf} and  also Haas~\cite{h:rfdsf} when  $\alpha'$ is
negative.

The  goal  of  this  paper  is   to  study  the  same  problem  for  the
fragmentation at nodes of the  L\'evy continuous random tree 
constructed in \cite{ad:falp}.

In  \cite{lglj:bplpep}  and \cite{lglj:bplplfss},  Le  Gall  and Le  Jan
associated to a Lévy process with  no negative jumps that does not drift
to  infinity,  $X=(X_s,  s\geq  0)$  with  Laplace  exponent  $\psi$,  a
continuous state  branching process (CSBP) and a  Lévy continuous random
tree (CRT) which keeps track of the genealogy of the CSBP.  The Lévy CRT
can  be coded  by  the so  called  height process,  $H=(H_s, s\geq  0)$.
Informally  $H_s$ gives  the distance  (which can  be understood  as the
number of generations) between the  individual labeled $s$ and the root,
0, of the CRT.  The precise definition of $\psi$ we consider is given at
the beginning of Section \ref{sec:exp}.

The  ideas  of \cite{ad:falp}  in  order  to  construct a  fragmentation
process from  this CRT  is to mark the  nodes of  the tree in  a Poissonian
manner. We then  cut the CRT at these marked nodes  and the ``sizes'' of
the  resulting subtrees  give the  state  of the  fragmentation at  some
time. As time $\theta$ increases, the parameter of the Poisson processes
used  to mark  the nodes  increases as  well as  the set  of  the marked
nodes. This gives a fragmentation process with no loss of mass.
When the initial L\'evy process is stable i.e. when  $\psi(\lambda)=\lambda^\alpha$,  $\alpha\in (1,2]$,
the fragmentation  is self-similar with  index $1/\alpha$ and  with a
zero  erosion  coefficient, see also see \cite{ap:sac} and
\cite{b:fpcbm} for $\alpha=2$, or \cite{m:sfdfstsn} for $\alpha \in
(1,2)$.
For a  general sub-critical  or critical CRT,  there is no  more scaling
property, and the dislocation measure, which describes how a fragment of
size  $r>0$ is  cut in  smaller  pieces, cannot  be expressed  as a  nice
function  of  the  dislocation measure  of  a  fragment  of size  1.  In
\cite{ad:falp},  the  authors give  the  family  of dislocation  measures
$(\nu_r, r>0)$ for the fragmentation at node of a general sub-critical or
critical    CRT.     Intuitively $\nu_r$ describes the way a mass $r$
breaks in smaller pieces.

We denote by $\N$ the excursion measure of the L\'evy process
$X$ (the fragmentation process is then defined under this measure). We
denote by $\sigma$ the length of the excursion. 
We have   (see Section
3.2.2. in \cite{dlg:rtlpsbp}) that
\begin{equation}
   \label{eq:exp-sigma}
\N[1-\expp{-\lambda \sigma}]=\psi^{-1}(\lambda),
\end{equation}
and $\psi^{-1}$ is the Laplace exponent of a subordinator (see
\cite{b:pl},  chap. VII), whose Lévy measure we denote by $\pi_*$. The
distribution of  $\sigma$  under $\N$ is given by $\pi_*$. As $\pi_*$ is
a Lévy measure, we have $\int_{(0,\infty )} (1\wedge r) \;
\pi_*(dr)<\infty $. 
For $\varepsilon>0$, we write 
\[
\bar \pi_*(\varepsilon)=\pi_*((\varepsilon,\infty ))=
\N[\sigma>\varepsilon]\quad\text{and}\quad
\varphi(\varepsilon)=\int_{(0,\varepsilon]} 
r\pi_*(dr)=\N[\sigma\ind_{\{\sigma\leq \varepsilon\}}]. 
\]

 If
$\Lambda(\theta)=(\Lambda_1(\theta),\Lambda_2(\theta),\ldots)$   is  the
state   of  the   fragmentation   at  time   $\theta$,   we  denote   by
$N^{\varepsilon}(\theta)$  the number  of fragments  of size  greater than
$\varepsilon$ i.e.
\[
N^{\varepsilon}(\theta)=\sum_{k=1}^{+\infty}\ind_{\{\Lambda_k(\theta)>\varepsilon\}}=\sup\{k\ge
1,   \Lambda_k(\theta)>\varepsilon\}
\]
 with   the   convention   $\sup
\emptyset=0$. And we  denote by $M^{\varepsilon}(\theta)$ the mass of
the fragments of size less than $\varepsilon$ i.e.
\[
M^{\varepsilon}(\theta)=\sum_{k=1}^{+\infty}\Lambda_k(\theta)\ind_{\{\Lambda_k(\theta)\le\varepsilon\}}=\sum_{k=N^{\varepsilon}(\theta)+1}^{+\infty}\Lambda_k(\theta). 
\]

Let $\cj=\{s\ge 0,\ X_s>X_{s-}\}$ and let $(\Delta_s,\ s\in\cj)$ be
the set of jumps of $X$. Conditionally on $(\Delta_s,s\in\cj)$, let
$(T_s,s\in\cj)$ be a family of independant random variables, such that
$T_s$ has exponential distribution with mean $1/\Delta_s$. $T_s$ is
the time at which the node of the CRT associated
to the jump $\Delta_s$ is marked in order to construct the
fragmentation process. 
Under $\N$, we denote by $R(\theta)$
the mass of the marked  nodes of the L\'evy CRT i.e.
$$R(\theta)=\sum_{s\in\cj\cap [0,\sigma]}\Delta_s\ind_{\{T_s\le \theta\}}.$$

The  main  result  of
this  paper  is then  the  following  Theorem. 

\begin{theo} 
\label{theo:smallN}
We have  $\displaystyle \lim_{\varepsilon\rightarrow 0}
\frac{N^{\varepsilon}(\theta)}{\bar \pi_*(\varepsilon)} =\lim_{\varepsilon
  \rightarrow 0} \frac{M^{\varepsilon}(\theta)}{\varphi(\varepsilon) }
=R(\theta)$ in $L^2\bigl(\N[\expp{-\beta\sigma}\cdot]\bigr)$, for any
$\beta>0$.
\end{theo}

We  consider  the   stable  case  $\psi(\lambda)=\lambda^\alpha$,  where
$\alpha\in      (1,2)$.       We      have      
$$\pi_*(dr)=      (\alpha
\Gamma(1-\alpha^{-1}))^{-1}     r^{-1-1/\alpha}     \;    dr,$$ 
which gives    
$$\bar
\pi_*(\varepsilon)= \Gamma(1-\alpha^{-1})^{-1} \varepsilon^{-1/\alpha},\qquad
\mbox{and} \qquad              
\varphi(\varepsilon)=                \Big((\alpha-1)
\Gamma(1-\alpha^{-1})\Big)^{-1}  \varepsilon^{1-\alpha^{-1}}.$$ 
{F}rom
scaling property, there exists
a version of $(\N_r, r>0)$ such that for all $r>0$ we have $\N_r[F((X_s,
s\in  [0,r]))]=\N_1[F((r^{1/\alpha}  X_{s/r} ,  s\in  [0,r]))]$ for  any
non-negative  measurable function  $F$  defined on  the  set of  càd-làg
paths.

\begin{prop}\label{prop:cv-a-stable}
Let $\psi(\lambda)=\lambda^\alpha$, for  $\alpha\in (1,2)$. For all
$\theta>0$, we have $\N$.a.e or $\N_1$-a.s.
\begin{equation}
   \label{eq:cv-stable}
\lim_{\varepsilon\to
  0}\Gamma(1-1/\alpha)\varepsilon^{1/\alpha}N^{\varepsilon}(\theta)=\lim_{\varepsilon\to
  0}(\alpha-1)\Gamma(1-1/\alpha)\frac{M^{\varepsilon}(\theta)}{\varepsilon^{1-1/\alpha}}=R(\theta).
\end{equation}
\end{prop}

\begin{rem}
  Notice  the   similarity  with   the  results  in   \cite{d:fhalp}  on
  asymptotics for  the small fragments  in case of the  fragmentation at
  height  of the  CRT: the  local  time of  the height  process is  here
  replaced by the functional $R$.
\end{rem}

\begin{rem}
\label{rem:b}
Let us compare the result of Proposition \ref{prop:cv-a-stable} with the
main Theorem of  \cite{b:smssf}, which we recall now. Let 
 $\Lambda$ be a self-similar fragmentation of index $\alpha>0$,
erosion coefficient $c=0$ and dislocation measure $\nu$. We set
\begin{align*}
\varphi_b(\varepsilon) & =\int_{\cs^\downarrow} (\sum_{i=1}^\infty  
\ind_{\{x_i>\varepsilon\}} -1) \nu_1(dx),\\
   f_b(\varepsilon)&=\int_{\cs^\downarrow} \sum_{i=1}^\infty  x_i
\ind_{\{x_i<\varepsilon\}}  \nu_1(dx),\\
g_b(\varepsilon) &=\int_{\cs^\downarrow} \left(\sum_{i=1}^\infty  x_i
\ind_{\{x_i<\varepsilon\}} \right)^2 \nu_1(dx).
\end{align*}
If there exists $\beta\in(0,1)$ such that $\varphi_b$ is regularly
varying at 0 with index $-\beta$ (which is equivalent to $f_b$ is
regularly varying at 0 with index $1-\beta$), and if there exists two positive
constants $c,\eta$ such that
\begin{equation}\label{eq:cond_bertoin}
g_b(\varepsilon)\le cf_b^2(\varepsilon)(\log
1/\varepsilon)^{-(1+\eta)},
\end{equation}
then a.s. 
$$\lim_{\varepsilon\to0}\frac{N^{\varepsilon}(\theta)}{\varphi_b(\varepsilon)}=\lim_{\varepsilon\to0}\frac{M^{\varepsilon}(\theta)}{f_b(\varepsilon)}=\int_0^\theta\sum_{i=1}^\infty
\Lambda_i(u)^{\alpha+\beta}du.$$

In our case, we have $\varphi$ and $\bar\pi_*$ equivalent to $\varphi_b$
and $f_b$ (up to  multiplicative constants, see Lemmas \ref{lem:f_b} and
\ref{lem:phi_b}). The normalizations are consequently the same. However,
we   have  here   $g_b(\varepsilon)=O(f_b^2(\varepsilon))$   (see  Lemma
\ref{lem:gOf2}) and  Bertoin's assumption \reff{eq:cond_bertoin}  is not
fulfilled. When this last assumption  holds, remark the limit process is
an increasing continuous process (as  $\theta$ varies). In our case this
assumption does not hold and the limit process $(R(\theta),\theta\ge 0)$
is  still increasing  but discontinuous  as $R(\theta)$  is a  pure jump
process  (this is  an increasing  sum of  marked masses).
\end{rem}

The paper is  organized as follows. In Section  \ref{sec:not}, we recall
the definition  and properties of  the height and  exploration processes
that  code  the  L\'evy  CRT  and  we recall  the  construction  of  the
fragmentation  process associated  to the  CRT.  The  proofs  of Theorem
\ref{theo:smallN}  and Proposition  \ref{prop:cv-a-stable} are  given in
Section  \ref{sec:smallN}. Notice  computations  given in  the proof  of
Lemma \ref{lem:NMLR-cond} based  on Propositions \ref{prop:cond-F-R} and
\ref{prop:cond-F} are  enough to  characterize the transition  kernel of
the  fragmentation $\Lambda$.  We  characterize the  law of  the scaling
limit $R(\theta)$ in Section \ref{sec:lawR}.  The computation needed for
Remark \ref{rem:b} are given in Section \ref{sec:appendix}.

\section{Notations}

 \label{sec:not}

\subsection{The exploration process}
\label{sec:exp}
Let $\psi$ denote the Laplace
exponent            of            $X$:           $\E\left[\expp{-\lambda
    X_t}\right]=\expp{t\psi(\lambda)}$,  $\lambda>0$.   We shall  assume
there is no Brownian part, so that
\[
\psi(\lambda)=\alpha_0\lambda+\int_{(0,+\infty)}\pi(d\ell)
\left[\expp{-\lambda\ell}-1+\lambda\ell\right],  
\]
with  $\alpha_0\ge  0$  and  the   Lévy  measure  $\pi$  is  a  positive
$\sigma$-finite measure  on $(0,+\infty)$ such  that $\int_{(0,+\infty)}
(\ell\wedge \ell^2)\pi(d\ell)<\infty$.  Following \cite{dlg:rtlpsbp}, we
shall also assume that $X$  is of infinite variation a.s.  which implies
that  $\int_{(0,1)}\ell\pi(d\ell)=\infty$.  Notice those  hypothesis are
fulfilled in the stable case: $\psi(\lambda)=\lambda^\alpha$, $\alpha\in
(1,2)$.      For      $\lambda\geq     1/\varepsilon>0$,     we     have
$\expp{-\lambda\ell}-1+ \lambda\ell\geq \frac{1}{2}\lambda \ell \ind_{\{
  \ell\geq   2\varepsilon\}}$,    which   implies   that   $\lambda^{-1}
\psi(\lambda)  \geq  \alpha_0+  \int_{(2\varepsilon,\infty  )}  \ell  \;
\pi(d\ell)$. We deduce that
\begin{equation}
   \label{eq:psi/l}
\lim_{\lambda
\rightarrow\infty } \frac{\lambda}{\psi(\lambda)} =0.
\end{equation}

The   so-called  exploration  process   $\rho=(\rho_t,t\ge  0)$   is
Markov process taking values in $\cm_f$, the set of positive measures on
$\R_+$.  The height process  at time $t$
is defined as  the supremum of the closed support  of $\rho_t$ (with the
convention  that  $H_t=0$ if  $\rho_t=0$).   Informally, $H_t$  gives the  distance (which  can be
understood as the number  of generations) between the individual labeled
$t$ and the root, 0, of the CRT.  In some sense $\rho_t(dv)$ records the
``number'' of brothers, with labels  larger than $t$, of the ancestor of
$t$ at generation $v$.

We recall the definition and properties of the exploration process which
are    given    in    \cite{lglj:bplpep},    \cite{lglj:bplplfss}    and
\cite{dlg:rtlpsbp}. The results of  this section are mainly extracted from
\cite{dlg:rtlpsbp}.
 
Let $I=(I_t,t\ge 0)$ be the infimum process of $X$, $I_t=\inf_{0\le s\le
  t}X_s$.  We will also consider for  every $0\le s\le t$ the infimum of
$X$ over $[s,t]$:
\[
I_t^s=\inf_{s\le r\le t}X_r.
\]
There exists a sequence $(\varepsilon_n,n\in \N^*)$ of positive real
numbers decreasing to 0 s.t. 
\[
\tilde H_t= \lim_{k\rightarrow\infty } \inv{\varepsilon_k} \int_0^t
\ind_{\{X_s<I^s_t+\varepsilon_k\}}\; ds
\]
exists and is finite a.s. for all $t\geq 0$. 

The point 0  is regular for the Markov process $X-I$,  $-I$ is the local
time  of $X-I$  at  0 and  the right  continuous  inverse of  $-I$ is  a
subordinator  with   Laplace  exponent  $\psi^{-1}$   (see  \cite{b:pl},
chap.  VII).    Notice  this  subordinator   has  no  drift   thanks  to
\reff{eq:psi/l}.  Let $\pi_*$ denote the corresponding Lévy measure.
Let $\N$ be the associated excursion measure of the process $X-I$ out of
0, and $\sigma=\inf\{t>0; X_t-I_t=0\}$ be the length of the excursion of
$X-I$ under  $\N$. Under  $\N$, $X_0=I_0=0$.

For $\mu\in \cm_f$, we define  $H^\mu=\sup\{x\in \supp \mu\}$, where
$\supp\mu$ is the 
closed   support  of   the  measure   $\mu$.  {F}rom   Section   1.2  in
\cite{dlg:rtlpsbp},   there    exists   a   $\cm_f$-valued   process,
$\rho^0=(\rho^0_t, t\geq 0)$, called  the exploration process, such that:
\begin{itemize}
\item A.s., for every $t\geq 0$, we have 
$\langle \rho_t^0,1\rangle =X_t-I_t$, and the process $\rho^0$ is càd-làg. 
\item The process $(H_s^0=H^{\rho^0_s}, s\geq 0)$ taking values in
  $[0,\infty ]$ is lower semi-continuous. 
\item   For  each   $t\geq   0$,  a.s.   $H^0_t=   \tilde  H_t$.
\item For every measurable non-negative function $f$ defined on $\R_+ $,
\[
\langle \rho^0_t,f\rangle =\int_{[0,t]} f(H^0_s)\; d_sI_t^s,
\]
or equivalently, with $\delta_x$ being the Dirac mass at $x$, 
\[
\rho^0_t(dr)=\sum_{\stackrel{0<s\le t}
  {X_{s-}<I_t^s}}(I_t^s-X_{s-})\delta_{H_s^0}(dr). 
\]
\end{itemize}

In the definition  of the exploration process, as $X$  starts from 0, we
have  $\rho_0=0$ a.s.  To  get the  Markov property  of $\rho$,  we must
define  the  process  $\rho$  started  at any  initial  measure  $\mu\in
\cm_f$.  For  $a\in [0,  \langle \mu,1\rangle ]  $, we  define the
erased measure $k_a\mu$ by
\[
k_a\mu([0,r])=\mu([0,r])\wedge (\langle \mu,1\rangle -a), \quad
\text{for $r\geq 0$}. 
\]

If $a> \langle  \mu,1\rangle $, we set $k_a\mu=0$.   In other words, the
measure $k_a\mu$ is the measure $\mu$ erased by a mass $a$ backward from
$H^\mu$.

For $\nu,\mu \in \cm_f $, and $\mu$ with compact support, we define
the concatenation $[\mu,\nu]\in \cm_f  $ of the two measures by:
\[
\bigl\langle [\mu,\nu],f\bigr\rangle =\bigl\langle \mu,f\bigr\rangle
+\bigl\langle \nu,f(H^\mu+\cdot)\bigr\rangle , 
\]
for $f$ non-negative measurable. Eventually, we  set for every $\mu\in \cm_f $ 
and every $t>0$, 
\[
\rho_t=\bigl[k_{-I_t}\mu,\rho_t^0].
\]
We say  that $\rho=(\rho_t, t\geq 0)$  is the process  $\rho$ started at
$\rho_0=\mu$, and  write $\P_\mu$ for its law.  We set $H_t=H^{\rho_t}$.
The process $\rho$ is càd-làg (with respect to the weak convergence
topology on $\cm_f$) and strong Markov.

\subsection{Notations for the fragmentation at nodes}
We recall the construction of the fragmentation under $\N$ given in \cite{ad:falp}
in an equivalent but easier way to understand. Recall $(\Delta_s, s\in \cj)$ is the set of jumps of
$X$ and $T_s$ is the time at which the jump $\Delta_s$ is marked. Conditionally on   $(\Delta_s, s\in \cj)$, $(T_s, s\in \cj)$ is a
family of independent random variables, such that $T_s$ has
exponential distribution
with mean $1/\Delta_s$. We consider the family of  measures (increasing in
$\theta$) defined for $\theta\geq 0$ and $t\geq 0$ by 
\[
\tilde m_t^\theta(dr)= 
\sum_{\stackrel{0<s\le t}{X_{s-}<I_t^s}} \; 
\ind_{\{T_s\leq \theta\}}  \; \delta_{H_s}(dr).
\]
Intuitively,  $\tilde m_t^\theta$  describes  the marked  masses of  the
measure  $\rho_t$ i.e.   the marked  nodes of  the associated  CRT.  

Then we cut  the CRT according  to these marks to  obtain the
state of the fragmentation process at time $\theta$.
To construct the fragmentation,
let  us  consider  the
following  equivalence  relation  $\mathcal{R}^\theta$ on  $[0,\sigma]$,
defined under $\N$ or $\N_\sigma$  by
\begin{equation}
   \label{eq:def-R}
s\mathcal{R}^\theta t\iff
\tilde m_s^\theta\bigl([H_{s,t},H_s]\bigr)=\tilde m_t^\theta
\bigl([H_{s,t},H_t]\bigr)=0,
\end{equation}
where  $\displaystyle  H_{s,t}=\inf_{u\in[s,t]}H_u$. Intuitively, two points $s$ and $t$ belongs to the
same class of equivalence (i.e. the same fragment) at time $\theta$ if there is no cut on their lineage down
to their most recent common ancestor, that is
$\tilde m_s^\theta$ put no mass on $[H_{s,t}, H_s]$ nor 
$\tilde m_t^\theta$ on $[H_{s,t}, H_t]$. 
Notice cutting occurs on 
branching points, that is at node of the CRT. Each node of the CRT
correspond to a jump of the underlying Lévy process $X$. 
The fragmentation process at time $\theta$ is then the Lebesgue
measures (ranked in non-increasing order) of the equivalent classes of
$\mathcal{R}^\theta$. 

\begin{rem}
  In  \cite{ad:falp}, see  definition  (14), we  use  another family  of
  measures  $m^{(\theta)}_t$.   From  their  construction,  notice  that
  $\tilde m^\theta_t$  is absolutely continuous  w.r.t. $m^{(\theta)}_t$
  and  $ m^{(\theta)}_t$  is absolutely continuous
  w.r.t. $\tilde m^\theta_t$, if  we  take  $T_s=\inf  \{V_{s,u}, u>
  0\}$,  where 
  $\sum_{u>0} \delta_{V_{s,u}}  $ is a  Poisson point measure  on $\R_+$
  with   intensity   $\Delta_s\ind_{\{u>0\}}$,   see  Section   3.1   in
  \cite{ad:falp}. In particular $\tilde m^\theta_t$ and $m^{(\theta)}_t$
  define  the   same  equivalence   relation  and  therefore   the  same
  fragmentation.
\end{rem}

In order to index the fragments, we define the ``generation'' of a
fragment. For any $s\le  \sigma$, let  us define $H_s^0=0$  and
recursively  for $k\in 
\N$,
\[
H_s^{k+1}=\inf\Bigl\{u\ge 0,  \tilde m_s^\theta\bigl((H_s^k,u]\bigr)>0\Bigr\},
\]
with the usual convention $\inf\emptyset=+\infty$. We  set the
``generation'' of $s$ as 
 \[
K_s=\sup\{j\in\N,\ H_s^j<+\infty\}.
\]
Notice that if $s\mathcal{R}^\theta t$, then $K_s=K_t$. In particular
all elements of a fragment have the same ``generation''. We also call
this ``generation'' the ``generation'' of the fragment. 
Let $(\sigma^{i,k}(\theta), i\in I_k)$ be  the family of lengths of fragments in
``generation''  $k$.   Notice that $I_0$ is reduced to one point, say 0,
and we write 
\[
\tilde \sigma(\theta)=\sigma^{0,0}(\theta)
\]
for the fragment which contains the root. The joint law of $(\tilde
\sigma(\theta), \sigma)$ is given in Proposition 7.3 in \cite{ad:falp}.

Let
$(r^{j,k+1}(\theta), j\in J_{k+1})$  be the family of sizes  of the marked nodes
attached to the snake of ``generation'' $k$. More precisely,
$$\bigl\{r^{j,k+1}(\theta), j\in J_{k+1}\bigr\}=\bigl\{\Delta_s,
T_s\le \theta\quad\mbox{and}\quad K_s=k+1\bigr\}.$$ 
We set, for $k\in \N$,
\[
L_k(\theta)  =\sum_{i\in I_k} \sigma^{i,k}(\theta),\quad
N_k^\varepsilon(\theta)  =\sum_{i\in 
  I_k} \ind_{\{\sigma^{i,k}(\theta)>\varepsilon\}},\quad 
M_k^\varepsilon(\theta)=\sum_{i\in I_k} \sigma^{i,k}(\theta)
\ind_{\{\sigma^{i,k}\leq 
  \varepsilon\}},
\]
 and  we set,   for $k\in \N^*$,
$$R_k(\theta)=\sum_{j\in J_k} r^{i,k}(\theta).$$
We set $R_0=0$. 
Let us remark that we have  $\sigma=\sum_{k\geq 0} L_k(\theta)$, $N^{\varepsilon}(\theta)=\sum_{k\geq 0}
N_k^\varepsilon(\theta)$, $M^{\varepsilon}(\theta)=\sum_{k\ge
  0}M_k^\varepsilon(\theta)$ and $R(\theta)=R_k(\theta)$.

Let $\cf_k$ be the $\sigma$-field generated by
$((\sigma^{i,l}(\theta), i\in I_l), R_l(\theta))_{0\le l\le k}$. 
As a consequence of the special Markov property (Theorem 5.2 of
\cite{ad:falp}) and using the recursive construction of Lemma 8.6 of \cite{ad:falp}, we
have the following Propositions. 

\begin{prop}
   \label{prop:cond-F-R}
   Under $\N$, conditionally   on   $\cf_{k-1}$    and   $R_k(\theta)$,   $\sum_{i\in   I_k}
   \delta_{\sigma^{i,k}(\theta)}$  is distributed   as a  Poisson  point process
   with intensity $R_k(\theta) \N[d\tilde \sigma(\theta) ]$.
\end{prop}

\begin{prop}
   \label{prop:cond-F}
   Under $\N$, conditionally   on   $\cf_{k-1}$,   $\sum_{j\in   J_k}
   \delta_{r^{j,k}(\theta)}$  is distributed  as a  Poisson  point process
   with intensity $L_{k-1}(\theta) (1-\expp{-\theta r})\; \pi(dr)$. 
\end{prop}

\begin{rem}
\label{rem:law-L}
   Those Propositions allow to compute the  law of the
   fragmentation $\Lambda(\theta)$ for a given $\theta$ (see
   computations of Laplace transform in the proof of Lemma
   \ref{lem:NMLR-cond}.
\end{rem}

Let us recall that the key object in \cite{ad:falp} is the tagged 
fragment which contains the root. Recall its size is denoted by  $\tilde
\sigma(\theta)$. This fragment   corresponds to the
subtree of the initial CRT (after pruning) that contains the root. This
subtree  is
 a L\'evy CRT and the Laplace exponent of the associated L\'evy
process is 
\[
\psi_\theta(\lambda):=\psi(\lambda+\theta)-\psi(\theta),\quad\lambda\geq
0. 
\]
This implies  $\psi_\theta^{-1}(v)=\psi^{-1}(v+\psi(\theta)) - \theta$
and  we deduce from \reff{eq:psi/l}
\begin{equation}
   \label{eq:cv-pq}
\lim_{\lambda\rightarrow\infty } \psi^{-1}_\theta(\lambda)/ \lambda=0.
\end{equation}
We also have (see (3) for the first equality with $\psi$ replaced by
$\psi_\theta$) 
\begin{equation}
   \label{eq:ets}
\N\left[1 -\expp{-\beta\tilde \sigma(\theta)}
\right]=
\psi_\theta^{-1}(\beta)\quad\text{and}\quad\N\left[\tilde\sigma(\theta) 
         \expp{-\beta\tilde           \sigma(\theta)} 
       \right]=\inv{\psi_\theta'(\psi^{-1}_\theta(\beta))}. 
\end{equation}

\section{Proofs}
\label{sec:smallN}
We fix $\theta>0$. As $\theta$ is fixed, we will omit to mention the
dependence w.r.t. $\theta$ of the different quantities in this section:
for example we write $\tilde \sigma$ and $N^\varepsilon$ for $\tilde
\sigma(\theta)$ and $N^\varepsilon(\theta)$. We set
\[
\cn^\varepsilon=N^{\varepsilon}-\ind_{\{\tilde\sigma>\varepsilon\}}
\qquad\mbox{and}\qquad \cm^\varepsilon=M^{\varepsilon}-\tilde
    \sigma\ind_{\{\tilde \sigma\leq \varepsilon\}}.
\]

\subsection{Proof of Theorem \ref{theo:smallN}}
The poof  is in  four steps. In  the first  step we compute  the Laplace
transform of $(\cn^\varepsilon,\cm^\varepsilon,  R, \sigma)$. From there
we  could prove  the  convergence of  Theorem  \ref{theo:smallN} with  a
convergence  in probability  instead  of  in $L^2$.  However  we need  a
convergence  speed to get  the a.s.  convergence in  the $\alpha$-stable
case of Proposition \ref{prop:cv-a-stable}. In the second step, we check
the computed Laplace transform has  the necessary regularity in order to
derive    in     the    third     step    the    second     moment    of
$(\cn^\varepsilon,\cm^\varepsilon,  R)$  under $\N[\expp{-\beta  \sigma}
\cdot]$.  In the last  step we  check the  convergence statement  of the
second moment.

In a \textbf{first step}, we give the joint law under $\N$ of
$(\cn^\varepsilon,\cm^\varepsilon, R, \sigma)$ by computing for $x>0$,
$y>0$, $\beta>0$, $\gamma>0$,
\[
\N\left[\expp{-(x\cn^\varepsilon+y\cm^\varepsilon+ \gamma R+\beta\sigma)}\Big | \tilde
  \sigma\right]. 
\]
By monotone convergence, we have
\begin{equation}
   \label{eq:monot-MNR}
\N\left[\expp{-(x\cn^\varepsilon+y\cm^\varepsilon+ \gamma
    R+\beta\sigma)}\Big|\tilde \sigma
\right]=\lim_{n\to\infty}\N\left[\expp{-\left(\beta \tilde \sigma +
      \sum_{l=1}^n(xN_l^\varepsilon+yM_l^\varepsilon+\gamma 
      R_l+\beta L_l)\right)}\Big|\tilde \sigma\right].
\end{equation}
We define the function $H_{(x,y,\gamma)}$ by
\[
H_{(x,y,\gamma)}(c)=G\left(\gamma+
    \N\left[1-\expp{-(x\ind_{\{\tilde\sigma>\varepsilon\}}+ 
 y\tilde \sigma\ind_{\{\tilde\sigma\leq \varepsilon\}}
+c\tilde \sigma)} \right]\right),
\]
where for $a\geq 0$, 
\begin{equation}
   \label{eq:def-G}
G(a)= \int\pi(dr)\left(1-\expp{-\theta
  r}\right)\left(1-\expp{-a r}\right)=
\psi(\theta+a) -\psi(a)-\psi(\theta)=\psi_\theta(a) - \psi(a).
\end{equation}
Recall $\cf_k$ is the $\sigma$-field generated by
$((\sigma^{i,l}, i\in I_l), R_l)_{0\le l\le k}$. 
We then have the following Lemma.
\begin{lem}
\label{lem:NMLR-cond}
   For $x,y,\gamma\in \R_+$, $\varepsilon>0$, we have for $k\in \N^*$, 
\[
\N \left[\expp{-(xN_k^\varepsilon+yM^\varepsilon_k + cL_k+\gamma
      R_k)}\Big|\cf_{k-1}\right]= \expp{-H_{(x,y,\gamma)}(c)
    L_{k-1}}.
\]
\end{lem}

\begin{proof}
As a consequence of Proposition \ref{prop:cond-F-R}, we have
\[
\N\left[\expp{-(xN_k^\varepsilon+yM^\varepsilon_k + cL_k+\gamma
      R_k)}\Big|\cf_{k-1}, R_k\right]= 
\expp{-R_k(\gamma+ \N[1-\exp(-(x\ind_{\{\tilde \sigma>\varepsilon\}} + y
\tilde \sigma\ind_{\{\tilde \sigma\leq \varepsilon\}} +c\tilde \sigma))])} .
\]
As a consequence of Proposition \ref{prop:cond-F}, we have
\[
\N\left[\expp{-R_k(\gamma+ \N[1-\exp(-(x\ind_{\{\tilde \sigma>\varepsilon\}} + y
\tilde \sigma\ind_{\{\tilde \sigma\leq \varepsilon\}} +c\tilde \sigma))])} 
\Big|\cf_{k-1}\right]= 
\expp{-H_{(x,y,\gamma)} (c)L_{k-1}}.
\]
\end{proof}
We define  the constants $c_{(k)}$ by induction:
\[
c_{(0)}=0\qquad\mbox{and}\qquad
c_{(k+1)}=H_{(x,y,\gamma)}(c_{(k)}+\beta).
\]
An immediate backward induction yields (recall
$L_0=\tilde \sigma$): for every integer $n\geq 1$, we have
\[
\N\left[\expp{-\left(\sum_{l=1}^n(xN_l^\varepsilon+yM^\varepsilon_l+\gamma
    R_l+\beta L_l)\right)}\Big|\tilde \sigma \right]=\expp{-c_{(n)}\tilde
\sigma }.
\]
Notice the function $G$ is of class $\cc^\infty $ on $(0,\infty )$,
concave increasing and
the function 
\[
c\mapsto  \N\left[1-\expp{-(x\ind_{\{\tilde\sigma>\varepsilon\}}+ 
 y\tilde \sigma\ind_{\{\tilde\sigma\leq \varepsilon\}}+
(\beta+c)\tilde \sigma)} \right]
\]
is  of class  $\cc^\infty $  on  $[0, \infty  )$ and is 
concave increasing. This implies that $H_{(x,y,\gamma)}$ is concave
increasing and of class  $\cc^\infty $. 
Notice that
\[
x\ind_{\{\tilde\sigma>\varepsilon\}} 
+ y\tilde \sigma\ind_{\{\tilde\sigma\leq \varepsilon\}}
+c\tilde \sigma
\leq (\frac{x}{\varepsilon} +y +c)\tilde \sigma.
\]
In particular,  we have $H_{(x,y,\gamma)}(c)  \leq G(\gamma+
\psi_\theta^{-1}(\frac{x}{\varepsilon} +y +c) )$.
As 
$\lim_{a\rightarrow   \infty   }   G'(a)=   0$, this   implies   that
$\lim_{a\rightarrow\infty          }          G(a)/a=0$.           Since
$\lim_{\lambda\rightarrow\infty } \psi^{-1}_\theta(\lambda)=\infty $, 
we deduce thanks to \reff{eq:cv-pq} that
\begin{equation}
   \label{eq:cv-H}
\lim_{c\rightarrow\infty } \frac{H_{(x,y,\gamma)}(c)}{c}= 0.
\end{equation}
For $\gamma>0$, notice $H_{(x,y,\gamma)}(0)>0$. As 
the function $H_{(x,y,\gamma)}$ is increasing and
continuous, we deduce the
sequence  $(c_{(n)}, n\geq  0)$  is increasing and converges
to   the unique root, say $c'$,  of
$c=H_{(x,y,\gamma)} (c+\beta)$. 
And we deduce from \reff{eq:monot-MNR} that 
\begin{equation}
   \label{eq:Nce}
\N\left[\expp{- (x\cn^\varepsilon+y\cm^\varepsilon+\gamma R+\beta\sigma)
}  \Big|\tilde \sigma \right]
= \expp{-(\beta+c') \tilde
    \sigma }.
\end{equation}

In a \textbf{second step}, we look at the dependency of
the root of $c=H_{(x,y,\gamma)} (c+\beta)$  in $(x,y,\gamma)$.

  Let  $\varepsilon,  x, y,  \beta,  \gamma\in  (0,\infty )$  be
  fixed.  There  exists $a>0$  small  enough  such  that for  all  $z\in
  (-a,a)$, we  have $z\gamma+\N[1-\expp{-\beta \tilde \sigma/2}] >0$ and
  for all  $\tilde \sigma \geq 
  0$,
\[
z(
x\ind_{\{\tilde \sigma>\varepsilon\}} +y\tilde \sigma \ind_{\{\tilde
  \sigma\leq \varepsilon\}} ) + \beta\tilde \sigma\geq \beta \tilde\sigma/2.
\]
We consider the function $J$ defined on $(-\beta/2,\infty )\times
(-a,a)$
by 
\[
J(c,z)=H_{zx,zy,z\gamma}(c+\beta) -c.
\]
{From} the regularity of $G$, we  deduce 
the function $J$ is  of class $\cc^\infty $ on $(-\beta/2,\infty )\times
(-a,a)$ and the function $c\mapsto J(c,z)$   is
 concave. Notice that $J(0,z)>0$ for all $z\in (-a,a)$. This and
\reff{eq:cv-H} implies that there exists a unique solution $c(z)$ to the
equation $J(c,z)=0$  and that $\displaystyle  \frac{\partial J}{\partial
  c} (c(z),z)<0$  for all $z\in  (-a,a)$. The implicit  function Theorem
implies  the function  $z\mapsto  c(z)$  is of  class  $\cc^\infty $  on
$(-a,a)$.      In      particular,      we      have      $\displaystyle
c(z)=c_0+zc_1+\frac{z^2}{2} c_2 + o(z^2)$.  We deduce from \reff{eq:Nce}
that for all $z\in [0,a)$,
\[
\N\left[\expp{- z(x\cn^\varepsilon+y\cm^\varepsilon+\gamma R)  -\beta\sigma
}  \Big|\tilde \sigma \right]
= \expp{-(\beta+c(z)) \tilde
    \sigma }= \expp{-(\beta+c_0+zc_1+\frac{z^2}{2} c_2 + o(z^2)) \tilde
    \sigma }.
\]

In   a   \textbf{third  step},   we   investigate   the  second   moment
$\N\left[(x\cn^\varepsilon+y\cm^\varepsilon+\gamma       R)^2\expp{-\beta
    \sigma} \right]$.
Standard results on Laplace transforms, implies the second moment is
finite and 
\begin{equation}
   \label{eq:NM-2g}
\N\left[(x\cn^\varepsilon+y\cm^\varepsilon+\gamma R)^2\expp{ -\beta\sigma
}  \Big|\tilde \sigma \right]
=\expp{-(\beta+c_0) \tilde
    \sigma }(c_1^2\tilde \sigma -  c_2 )\tilde \sigma.
\end{equation}
Next we compute $c_0$, $c_1$ and $c_2$.
By definition of $c(z)$, we have
\[
c_0+zc_1+\frac{z^2}{2} c_2 + o(z^2)= G\left(z\gamma+ 
\N\left[1-\expp{-z(x\ind_{\{\tilde\sigma>\varepsilon\}}+ 
 y\tilde \sigma\ind_{\{\tilde\sigma\leq \varepsilon\}}) -
(\beta+c_0+zc_1+\frac{z^2}{2} c_2 + o(z^2))\tilde \sigma}
\right]\right).
\]
We compute the expansion in $z$ of the right hand-side term of this
equality. We set 
\begin{align*}
   a_0&= \N\left[1-\expp{-(\beta+c_0)\tilde
       \sigma)}\right]=\psi_\theta^{-1}(\beta+ c_0) ,\\
   a_1&=\gamma+ \N\left[\expp{-(\beta+c_0)\tilde
       \sigma}(x\ind_{\{\tilde\sigma>\varepsilon\}}+  
 y\tilde \sigma\ind_{\{\tilde\sigma\leq \varepsilon\}} +c_1 \tilde
     \sigma) \right],\\
   a_2&=\N\left[\expp{-(\beta+c_0)\tilde
       \sigma}(c_2\tilde \sigma - (x\ind_{\{\tilde\sigma>\varepsilon\}}+  
 y\tilde \sigma\ind_{\{\tilde\sigma\leq \varepsilon\}} +c_1 \tilde
     \sigma)^2)  \right], 
\end{align*}
so that  standard results on Laplace transform yield
\[
c_0+zc_1+\frac{z^2}{2} c_2 + o(z^2)= G\left(a_0+za_1+\frac{z^2}{2} a_2+
  o(z^2)\right).
\]
We deduce that 
\begin{align}
\nonumber%\label{eq:c0}
   c_0&=G(a_0)=G\left(\N\left[1 -\expp{-(\beta+c_0)\tilde \sigma}
       \right]\right),\\ 
\label{eq:c1}
c_1& =a_1 G'(a_0),\\
\label{eq:c2}
c_2&= a_2 G'(a_0)+ a_1^2 G''(a_0)= a_2G'(a_0)+\frac{c_1^2 G''(a_0)}{G'(a_0)^2}.
\end{align}
Using \reff{eq:def-G} and \reff{eq:ets}, we have
$c_0=G\left(\psi_\theta^{-1}(\beta+c_0)\right)=\beta+c_0 - 
\psi(\psi_\theta^{-1}(\beta+c_0))$,
that is 
\[
h_\beta:=\beta+c_0=\psi_\theta(\psi^{-1}(\beta))\quad\text{and}\quad 
a_0=\psi_\theta^{-1}(\beta+ c_0)=\psi^{-1}(\beta). 
\]
Notice that $h_\beta>0$. 
And we have, thanks to the second equality of \reff{eq:ets}, 
\[
G'(\psi^{-1}(\beta)) \N\left[\expp{-h_\beta\tilde
       \sigma}\tilde \sigma  \right] = \frac{\psi_\theta'(\psi^{-1}(\beta))
   -\psi'(\psi^{-1}(\beta))}{\psi_\theta'(\psi^{-1}(\beta))}<1. 
\]
(This  last   inequality  is  equivalent  to   say  that  $\displaystyle
\frac{\partial   J}{\partial   c}    (c(z),z)<0$   at   $z=0$.)    From
\reff{eq:c1}, we get
\begin{equation}
   \label{eq:c1b}
c_1= G'(\psi^{-1}(\beta))\frac{\gamma+ \N\left[\expp{-h_\beta\tilde
       \sigma}(x\ind_{\{\tilde\sigma>\varepsilon\}}+  
 y\tilde \sigma\ind_{\{\tilde\sigma\leq \varepsilon\}}) \right]}
{1- G'(\psi^{-1}(\beta)) \N\left[\expp{-h_\beta\tilde
       \sigma}\tilde \sigma  \right]},
\end{equation}
and from \reff{eq:c2}
\begin{equation}
   \label{eq:c2b}
c_2=\frac{-G'(\psi^{-1}(\beta)) \N\left[\expp{-h_\beta\tilde
       \sigma}(x\ind_{\{\tilde\sigma>\varepsilon\}}+  
 y\tilde \sigma\ind_{\{\tilde\sigma\leq \varepsilon\}} +c_1 \tilde
     \sigma)^2 \right]  +  \frac{\displaystyle c_1^2
     G''(\psi^{-1}(\beta))}{\displaystyle G'(\psi^{-1}(\beta))^{2} }}{ 1- G'(\psi^{-1}(\beta)) \N\left[\expp{-h_\beta\tilde
       \sigma}\tilde \sigma  \right]}.
\end{equation}
We get 
\begin{equation}
   \label{eq:MN-2b}
\N\left[(x\cn^\varepsilon+y\cm^\varepsilon+\gamma R)^2\expp{-\beta
       \sigma} \right]= c_1^2\N\left[\expp{-h_\beta\tilde
       \sigma}\tilde \sigma^2  \right]  - c_2 \N\left[\expp{-h_\beta\tilde
       \sigma}\tilde \sigma  \right], 
\end{equation}
where  $c_1$ and $c_2$  defined by  \reff{eq:c1b} and  \reff{eq:c2b} are
polynomials of respective degree 1 and 2 in $x,y$ and $\gamma$. In particular
\reff{eq:MN-2b} also holds for $x,y,\gamma\in \R$. 

In a \textbf{fourth step}, we look at asymptotics as $\varepsilon$
decreases to 0. Let $\lambda_1,  \lambda_2\in \R_+$ and $\gamma=-(\lambda_1+\lambda_2)$. We  set 
\[
x_\varepsilon=\lambda_1/
  \bar                      \pi_*(\varepsilon)\quad \text{and} \quad 
y_\varepsilon=\lambda_2/\varphi(\varepsilon).
\]
We recall from Lemma 4.1 in  \cite{d:fhalp}, that 
\begin{equation}
   \label{eq:cv-f-e}
\lim_{\varepsilon \rightarrow 0} \inv{\bar \pi_*(\varepsilon)}=0
\quad\text{and} \quad \lim_{\varepsilon \rightarrow 0}
\frac{\varepsilon}{\varphi(\varepsilon)}
=0.
\end{equation}
Lemma 7.2 in \cite{ad:falp} tells that for any non-negative measurable
function $F$, we have 
$\N[F(\tilde \sigma)]=\N[\expp{-\theta \sigma}F(\sigma)]$. 
We define 
\begin{align*}
\Delta_\varepsilon
&:=\gamma+ \N\left[\expp{-h_\beta\tilde
       \sigma}(x_\varepsilon\ind_{\{\tilde\sigma>\varepsilon\}}+  
 y_\varepsilon\tilde \sigma\ind_{\{\tilde\sigma\leq \varepsilon\}})
\right]\\
&=  - x_\varepsilon \N\left[(1-\expp{-(h_\beta+\theta)
       \sigma})\ind_{\{\sigma>\varepsilon\}}\right]-  
 \varepsilon y_\varepsilon\N\left[(1-\expp{-(h_\beta+\theta)
       \sigma})\frac{ \sigma}{\varepsilon}
     \ind_{\{\sigma\leq \varepsilon\}}) \right].
\end{align*}
In particular, we have $\displaystyle \Delta_\varepsilon=O\left( \inv{\bar
  \pi_*(\varepsilon)} + \frac{\varepsilon}{\varphi(\varepsilon)}\right)$ and
from \reff{eq:cv-f-e} $\lim_{\varepsilon\rightarrow 0}
\Delta_\varepsilon=0$. 
{From} \reff{eq:c1b}, we get  $\displaystyle
c_1=O(\Delta_\varepsilon)= O\left( \inv{\bar 
  \pi_*(\varepsilon)} + \frac{\varepsilon}{\varphi(\varepsilon)}\right)$. 
{From} \reff{eq:c2b}, we also have for some finite constant $C$
independent of $\varepsilon$:
\[
|c_2|\leq 2\frac{G'(\psi^{-1}(\beta)) \N\left[\expp{-h_\beta\tilde
       \sigma}(x_\varepsilon\ind_{\{\tilde\sigma>\varepsilon\}}+  
 y_\varepsilon\tilde \sigma\ind_{\{\tilde\sigma\leq \varepsilon\}})^2 \right]  }{ 1- G'(\psi^{-1}(\beta)) \N\left[\expp{-h_\beta\tilde
       \sigma}\tilde \sigma  \right]}+ Cc_1^2. %O(\Delta_\varepsilon).
\]
Notice that 
\begin{align}
\nonumber
   \N\left[\expp{-h_\beta\tilde
       \sigma}(x_\varepsilon\ind_{\{\tilde\sigma>\varepsilon\}}+  
 y_\varepsilon\tilde \sigma\ind_{\{\tilde\sigma\leq \varepsilon\}})^2
\right]
&=\N\left[\expp{-h_\beta\tilde
       \sigma}(x_\varepsilon^2\ind_{\{\tilde\sigma>\varepsilon\}}+  
 y_\varepsilon^2\tilde \sigma^2\ind_{\{\tilde\sigma\leq \varepsilon\}})
\right]\\
\nonumber
&\leq \frac{\lambda_1}{\bar \pi_*(\varepsilon)}
\;x_\varepsilon\N\left[\expp{-h_\beta\tilde 
       \sigma}\ind_{\{\tilde\sigma>\varepsilon\}}\right]
+ \frac{\lambda_2\varepsilon}{\varphi(\varepsilon)} \; y_\varepsilon\N\left[\expp{-h_\beta\tilde
       \sigma}
\tilde \sigma \ind_{\{\tilde\sigma\leq \varepsilon\}}
\right]\\
&=O\left( \inv{\bar
  \pi_*(\varepsilon)} + \frac{\varepsilon}{\varphi(\varepsilon)}\right).
\label{eq:majoN}
\end{align}
We deduce that $\displaystyle c_2=O\left( \inv{\bar
  \pi_*(\varepsilon)} + \frac{\varepsilon}{\varphi(\varepsilon)}\right)$. 
Equation \reff{eq:MN-2b}  implies that 
\[
\N\left[\left(\lambda_1\frac{\cn^\varepsilon}{\bar \pi_*(\varepsilon)}
  +\lambda_2
  \frac{\cm^\varepsilon}{\varphi(\varepsilon)}-(\lambda_1+\lambda_2)
  R\right)^2\expp{-\beta 
       \sigma} \right]=O\left( \inv{\bar
  \pi_*(\varepsilon)} + \frac{\varepsilon}{\varphi(\varepsilon)}\right).
\]
As $\sigma\geq \tilde \sigma$, we have 
\begin{align*}
\N\left[\left(\inv{\bar \pi_*(\varepsilon)^2} \ind_{\{\tilde \sigma
      >\varepsilon\}} +  \inv{\varphi(\varepsilon)^2}\tilde
    \sigma^2\ind_{\{\tilde \sigma 
      \leq \varepsilon\}} \right)\expp{-\beta 
       \sigma} \right]
&\leq \N\left[\left(\inv{\bar \pi_*(\varepsilon)^2} \ind_{\{\tilde \sigma
      >\varepsilon\}} +  \inv{\varphi(\varepsilon)^2}\tilde
    \sigma^2\ind_{\{\tilde \sigma 
      \leq \varepsilon\}} \right)\expp{-\beta  \tilde 
       \sigma} \right]\\
&= O\left( \inv{\bar
  \pi_*(\varepsilon)} + \frac{\varepsilon}{\varphi(\varepsilon)}\right),
 \end{align*}
where we used \reff{eq:majoN} for the last equation (with $\beta$
instead of $h_\beta$).  Recall that $N^{\varepsilon}(\theta)=\cn^\varepsilon + \ind_{\{\tilde \sigma
  >\varepsilon\}} $ and $M^{\varepsilon}(\theta)= \cm^\varepsilon+ \tilde \sigma
\ind_{\{\tilde \sigma   \leq \varepsilon\}}$ and thus
\begin{multline*}
   \N\left[\left(\lambda_1\frac{N^\varepsilon}{\bar \pi_*(\varepsilon)}
  +\lambda_2
  \frac{M^\varepsilon}{\varphi(\varepsilon)}-(\lambda_1+\lambda_2)
  R\right)^2\expp{-\beta 
       \sigma} \right]\\
\leq 2 \N\left[\left(\lambda_1\frac{\cn^\varepsilon}{\bar \pi_*(\varepsilon)}
  +\lambda_2
  \frac{\cm^\varepsilon}{\varphi(\varepsilon)}-(\lambda_1+\lambda_2)
  R\right)^2\expp{-\beta 
       \sigma} \right] \\
+ 2(\lambda_1^2+\lambda_2^2)\N\left[\left(\inv{\bar \pi_*(\varepsilon)^2} \ind_{\{\tilde \sigma
      >\varepsilon\}} +  \inv{\varphi(\varepsilon)^2}\tilde
    \sigma^2\ind_{\{\tilde \sigma 
      \leq \varepsilon\}} \right)\expp{-\beta 
       \sigma} \right]
\end{multline*}  
We deduce that  
\begin{equation}
   \label{eq:cvNM2}
\N\left[\left(\lambda_1\frac{N^\varepsilon}{\bar \pi_*(\varepsilon)}
  +\lambda_2
  \frac{M^\varepsilon}{\varphi(\varepsilon)}-(\lambda_1+\lambda_2)
  R\right)^2\expp{-\beta 
       \sigma} \right]=O\left( \inv{\bar
  \pi_*(\varepsilon)} + \frac{\varepsilon}{\varphi(\varepsilon)}\right).
\end{equation}
which, thanks to \reff{eq:cv-f-e},  exactly says that 
$\displaystyle \lim_{\varepsilon\rightarrow 0}
\frac{N^{\varepsilon}(\theta)}{\bar \pi_*(\varepsilon)} =\lim_{\varepsilon
  \rightarrow 0} \frac{M^{\varepsilon}(\theta)}{\varphi(\varepsilon) } =R$ in
$L^2(\N[\expp{- \beta \sigma} \cdot])$. 

\subsection{Proof of Proposition \ref{prop:cv-a-stable}}
\label{sec:proof-prop}
Recall that, in the stable case, we have
\[
\bar\pi_*(\varepsilon)=\frac{1}{\Gamma(1-1/\alpha)}\varepsilon^{-1/\alpha}\qquad\mbox{and}\qquad\varphi(\varepsilon)=\frac{1}{(\alpha-1)\Gamma(1-1/\alpha)}\varepsilon^{1-1/\alpha}.
\]
    Therefore $\varepsilon_n=n^{-2\alpha}$, $n\geq 1$,
   satisfies $\displaystyle \sum_{n\geq 1} \inv{\bar\pi_*(\varepsilon_n)}+
   \frac{\varepsilon_n}{\varphi(\varepsilon_n)} <\infty 
   $. The series with general term given by the left hand-side of
   \reff{eq:cvNM2} with $\varepsilon=\varepsilon_n$ is convergent. This
   implies that $\N$-a.e (and $\N_1$-a.s.)
$$\lim_{n\rightarrow \infty }
\frac{N^{\varepsilon_n}(\theta)}{\bar \pi_*(\varepsilon_n)} =\lim_{n
  \rightarrow \infty } \frac{M^{\varepsilon_n}(\theta)}{\varphi(\varepsilon_n) } =R(\theta).$$
Since
   $N^{\varepsilon}(\theta)$ is a non-increasing function of $\varepsilon$, we
   get that for any $\varepsilon\in [(n+1)^{-2\alpha}, n^{-2\alpha}]$,
   we have 
\[
\frac{n^2}{(n+1)^2} n^{-2}
   N^{n^{-2\alpha}}(\theta) \leq  \varepsilon^{1/\alpha} N^{\varepsilon}(\theta)
\leq  \frac{(n+1)^2}{n^2} (n+1)^{-2}
   N^{(n+1)^{-2\alpha}}(\theta) .
\]
   Hence we deduce that  $\N$-a.e. or $\N_1$-a.s., $\lim_{\varepsilon
     \rightarrow\ 0 } \varepsilon^{1/\alpha} 
   N^{\varepsilon}(\theta) =R(\theta)/ \Gamma(1-\alpha^{-1})$. 

The proof for
   $M^{\varepsilon}(\theta)$ is similar, as $M^{\varepsilon}(\theta)$ is a
   non-decreasing function of $\varepsilon$.

\section{Law of $R(\theta)$}\label{sec:lawR}

\begin{lem}
   \label{lem:R_s=1}
Let $\beta\geq 0$, $\gamma\leq 0$. We have 
\[
\N\left[1- \expp{- \beta \sigma -\gamma R(\theta)
}  \right]
= v
\]
where $v$ is  the unique non-negative root of
\begin{equation}
   \label{eq:root-x}
\beta+\psi(\gamma+\theta+v)=\psi(v+\theta)+\psi(v+\gamma).
\end{equation}
\end{lem}

\begin{rem}
  For the limit case $\psi(\lambda)=\lambda^2$ (which is excluded here),
  we   get  the   unique  non-negative   root  of   \reff{eq:root-x}  is
  $v=\sqrt{\lambda+2\gamma\theta}$.    This    would   implies
  $R(\theta)=2\theta  \sigma$  $\N$-a.e.  and $R(\theta)=2\theta$  
  $\N_1$-a.s.  This agrees with  the result  in \cite{b:smssf},  where the
  limit  which   appears  for  \reff{eq:cv-stable}  is   a.s.  equal  to
  $2\theta$.
\end{rem}

\begin{proof}
Take $x=y=0$ in \reff{eq:Nce}, integrate w.r.t. $\N$ and use
\reff{eq:ets} to get  
\[
\N\left[1 -\expp{- \gamma R(\theta)-\beta\sigma
}  \right]
= \N\left[1 -\expp{-(\beta+c) \tilde
    \sigma(\theta) }\right] = \psi_\theta^{-1}(\beta+c),
\]
where $c$ is the unique root of $c=H_{(0,0,\gamma)}(c)$ that is of
$c=G(\gamma+ \psi_\theta^{-1}(\beta+c))$.  
If we set
$v=\psi_\theta^{-1} (\beta +c) $, we have 
that $v$ is  the unique non-negative root of the
equation $G(\gamma+v)=\psi_\theta(v)- \beta$, that is
\reff{eq:root-x}. 
\end{proof}

\section{Appendix}\label{sec:appendix}
Let $\alpha\in  (1,2)$.  
Recall  from
\cite{ad:falp} Corollary 9.3 or \cite{m:sfdfstsn} that the fragmentation
is self similar with index $1/\alpha$ and dislocation measure given by
\[
\int_{\cs^\downarrow} F(x) \nu_1(dx)= \frac{\alpha(\alpha-1)
  \Gamma(1-1/\alpha)}{\Gamma(2-\alpha)} \E[S_1 \;F((\Delta
S_t/ S_1, t\leq 1))],
\]
where $F$ is any non-negative   measurable    function  on
$\cs^\downarrow$, and   $(\Delta S_t,  t\geq 0)$ are  the jumps  of a
stable   subordinator   $S=(S_t,   t\geq   0)$   of   Laplace   exponent
$\psi^{-1}(\lambda)=\lambda^{1/\alpha}$, ranked by decreasing size.

In  this  section  we  shall  compute  the functions $f_b$, $\varphi_b$
and  $g_b$ defined in 
\cite{b:smssf}  and recalled in Remark \ref{rem:b} for the self-similar
fragmentation at nodes. 
\begin{lem}\label{lem:f_b}
We have
$\displaystyle f_b(\varepsilon)=\frac{1}{\Gamma(1+1/\alpha)}
\left(\frac{\varepsilon}{1-\varepsilon} \right)^{1-1/\alpha}.$
\end{lem}

\begin{proof}
The Lévy measure of $S$ is given by $\displaystyle
   \pi_*(dr)=\inv{\alpha    \Gamma    (1-1/\alpha)} \frac{dr}{r^{
1+1/\alpha}}  dr$. 
For $\beta\in (0,1)$, we have $\displaystyle \int_{(0,\infty )}
  \frac{dy}{y^{1+\beta} } (1- \expp{-y\lambda})=\lambda^\beta
  \frac{\Gamma(1-\beta)}{\beta}$.   We deduce that 
\[
   \E[S_1^\beta]
=\frac{\beta}{\Gamma(1-\beta)} \E\left[\int_0^\infty
  \frac{dy}{y^{1+\beta} } (1- \expp{-yS_1})\right]
= \frac{\beta}{\Gamma(1-\beta)} \int_0^\infty
  \frac{dy}{y^{1+\beta} } (1- \expp{-y^{1/\alpha}})
 = \frac{\Gamma(1-\alpha
    \beta)}{\Gamma(1-\beta)} . 
\]
  Standard computation for Poisson measure yield
\begin{align*}
   f_b(\varepsilon)
&=\int_{\cs^\downarrow} \sum_{i=1}^\infty  x_i
\ind_{\{x_i<\varepsilon\}}  \nu_1(dx)\\
&=  \frac{\alpha(\alpha-1)
  \Gamma(1-1/\alpha)}{\Gamma(2-\alpha)} \E\left[S_1 \sum_{t\leq 1} \frac{\Delta
S_t}{ S_1}\ind_{\{\Delta S_t <\varepsilon S_1\}}\right]\\ 
&=  \frac{\alpha(\alpha-1)
  \Gamma(1-1/\alpha)}{\Gamma(2-\alpha)} \E\left[ \sum_{t\leq 1} 
\Delta S_t\ind_{\{\Delta S_t <\varepsilon (S_1-\Delta S_t)/(1-\varepsilon)\}}\right]\\
&=  \frac{\alpha(\alpha-1)
  \Gamma(1-1/\alpha)}{\Gamma(2-\alpha)} \E\left[\int  \pi_*(dr) r 
\ind_{\{r <\varepsilon S_1/(1-\varepsilon)\}}\right]\\
&= \frac{\alpha}{\Gamma(2-\alpha)} \E[S_1^{1-1/\alpha}]
\left(\frac{\varepsilon}{1-\varepsilon} \right)^{1-1/\alpha}\\
&= \frac{1}{\Gamma(1+1/\alpha)}
\left(\frac{\varepsilon}{1-\varepsilon} \right)^{1-1/\alpha}. 
\end{align*}
\end{proof}

\begin{lem}\label{lem:phi_b}
We have
$\displaystyle \lim_{\varepsilon \rightarrow 0} \varepsilon^{1/\alpha} \varphi_b(\varepsilon)= \frac{\alpha-1
  }{\Gamma(1+1/\alpha)} $.
\end{lem}
\begin{proof}
We have 
\begin{align*}
 \varphi_b(\varepsilon)
& =\int_{\cs^\downarrow} (\sum_{i=1}^\infty  
\ind_{\{x_i>\varepsilon\}} -1) \nu_1(dx)\\
&=  \frac{\alpha(\alpha-1)
  \Gamma(1-1/\alpha)}{\Gamma(2-\alpha)} \E\left[S_1 \sum_{t\leq 1}
  \ind_{\{\Delta S_t > \varepsilon S_1\}}- S_1\right]\\ 
&=  \frac{\alpha(\alpha-1)
  \Gamma(1-1/\alpha)}{\Gamma(2-\alpha)} \E\left[ \sum_{t\leq 1}(S_1
  -\Delta S_t) 
  \ind_{\{\Delta S_t > \varepsilon \frac{S_1-\Delta S_t}{1-\varepsilon}\}}-
  \Delta S_t \ind_{\{\Delta S_t \leq  \varepsilon \frac{S_1-\Delta S_t}{
    1-\varepsilon}\}}\right]\\ 
&=\frac{\alpha(\alpha-1)
  \Gamma(1-1/\alpha)}{\Gamma(2-\alpha)}  \E\left[ S_1 \int \pi_*(dr)
  \ind_{\{ r> \varepsilon S_1/(1-\varepsilon)\}} \right] -
f_b(\varepsilon) \\ 
&=\frac{\alpha(\alpha-1)
  }{\Gamma(2-\alpha)}  \E\left[ S_1^{1-1/\alpha}
  \right]\left(\frac{\varepsilon}{1-\varepsilon} \right)^{-1/\alpha}-
  f_b(\varepsilon)\\  
&=\frac{\alpha-1
  }{\Gamma(1+1/\alpha)} \left(\frac{\varepsilon}{1-\varepsilon}
  \right)^{-1/\alpha}- f_b(\varepsilon)
\end{align*}
\end{proof}

\begin{lem}\label{lem:gOf2}
   The limit  $\displaystyle \lim_{\varepsilon \rightarrow
  0} \frac{g_b(\varepsilon)}{f_b(\varepsilon)^2}$ exists and belongs to
$(0,\infty )$.  
\end{lem}
\begin{proof}

We have 
\begin{align*}
   g_b(\varepsilon)
 &=\int_{\cs^\downarrow} \left(\sum_{i=1}^\infty  x_i
\ind_{\{x_i<\varepsilon\}} \right)^2 \nu_1(dx)\\
&=  \frac{\alpha(\alpha-1)
  \Gamma(1-1/\alpha)}{\Gamma(2-\alpha)} \E\left[S_1 \left(\sum_{t\leq 1} \frac{\Delta
S_t}{ S_1}\ind_{\{\Delta S_t <\varepsilon S_1\}}\right)^2 \right]\\ 
&=  \frac{\alpha(\alpha-1)
  \Gamma(1-1/\alpha)}{\Gamma(2-\alpha)} \Big(\E\Big[ \sum_{t\leq 1} 
\frac{(\Delta S_t)^2}{S_1} \ind_{\{\Delta S_t
  <\varepsilon S_1\}}\Big]  \\
 &\hspace{1cm}
+ \E\Big[ \sum_{t\leq 1, s\leq 1, s\neq t} 
\frac{\Delta S_t \Delta S_s }{S_1} \ind_{\{\Delta S_t <\varepsilon S_1, \Delta S_s <
  \varepsilon S_1\}}\Big] \Big). 
\end{align*}
For the first term, we get 
\begin{align*}
\E\Big[ \sum_{t\leq 1} 
\frac{(\Delta S_t)^2}{S_1} \ind_{\{\Delta S_t
  <\varepsilon S_1\}}\Big]
&\leq \E\Big[ \sum_{t\leq 1} 
\frac{(\Delta S_t)^2}{S_1-\Delta S_t} \ind_{\{\Delta S_t
  <\varepsilon (S_1-\Delta 
  S_t)/(1-\varepsilon)\}}\Big]\\
&=    \E\left[\inv{S_1} \int  \pi_*(dr) r^2  
\ind_{\{r <\varepsilon S_1/(1-\varepsilon)\}}\right]\\
&= \frac{1}{(2\alpha -1) \Gamma(1-1/\alpha)} \E[S_1^{1-1/\alpha}]
\left(\frac{\varepsilon}{1-\varepsilon} \right)^{2-1/\alpha}\\
&= o(\varepsilon^{2-2/\alpha}). 
\end{align*}
For the second term, we notice that for $r,s,S \in \R_+$
\begin{align*}
  \{ r\leq \varepsilon S/(1-\varepsilon), v\leq \varepsilon
S/(1-\varepsilon) \} 
&\subset 
\{ r\leq \varepsilon (S+r+v), v\leq \varepsilon (S+r+v) \}\\
&\subset 
\{ r\leq \varepsilon S/(1-2\varepsilon), v\leq \varepsilon
S/(1-2\varepsilon) \} .
\end{align*}
And we get 
\begin{multline*}
\E\Big[ \sum_{t\leq 1, s\leq 1, s\neq t} 
\frac{\Delta S_t \Delta S_s}{S_1} \ind_{\{\Delta S_t
  <\varepsilon S_1, \Delta S_s < \varepsilon S_1\}}\Big]\\
\begin{aligned}
&\leq \E\Big[ \sum_{t\leq 1, s\leq 1, s\neq t} 
\frac{\Delta S_t \Delta S_s}{S_1-\Delta S_t -\Delta S_s} \ind_{\{\Delta S_t
  <\varepsilon \frac {S_1- \Delta S_t -\Delta S_s}{1-2\varepsilon} ,  \Delta
  S_s < \varepsilon \frac{S_1- \Delta S_t -\Delta
  S_s}{1-2\varepsilon}\}}\Big] \\
&=\E\left[\inv{S_1} \left(\int \pi_*(dr)r \ind_{\{r<\varepsilon
      S_1/(1-2\varepsilon) \}} \right)^2\right]\\
&=c_\alpha \left(\frac{\varepsilon}{1-2\varepsilon} \right)^{2-2/\alpha},
\end{aligned}
\end{multline*}

with $\displaystyle c_\alpha=\frac{\Gamma(3-\alpha)}{(\alpha-1)^2\Gamma(2/\alpha)\Gamma(1-1/\alpha)^2}$,  as well as
\begin{multline*}
\E\Big[ \sum_{t\leq 1, s\leq 1, s\neq t} 
\frac{\Delta S_t \Delta S_s}{S_1} \ind_{\{\Delta S_t
  <\varepsilon S_1, \Delta S_s < \varepsilon S_1\}}\Big]\\
\begin{aligned}
&\geq \E\Big[ \sum_{t\leq 1, s\leq 1, s\neq t} 
\frac{\Delta S_t \Delta S_s}{(S_1-\Delta S_t -\Delta S_s)\frac{1+2\varepsilon}{1-\varepsilon}}
\ind_{\{\Delta S_t 
  <\varepsilon \frac{S_1- \Delta S_t -\Delta S_s}{1-\varepsilon} ,  \Delta
  S_s < \varepsilon \frac{S_1- \Delta S_t -\Delta
  S_s}{1-\varepsilon}\}}\Big] \\
&=\frac{1-\varepsilon}{1+2\varepsilon} \E\left[\inv{S_1} \left(\int \pi_*(dr)r \ind_{r<\varepsilon
      S_1/(1-\varepsilon) \}} \right)^2\right]\\
&=c_\alpha
\frac{1-\varepsilon}{1+2\varepsilon}\left(\frac{\varepsilon}{1-\varepsilon}
\right)^{2-2/\alpha}. 
\end{aligned}
\end{multline*}
In particular, we have that $g_b(\varepsilon)=c_\alpha
\varepsilon^{2-2/\alpha} (1+o(1))$. We deduce that  
$\displaystyle \lim_{\varepsilon \rightarrow
  0} \frac{g_b(\varepsilon)}{f_b(\varepsilon)^2}\in (0,\infty )$. 
\end{proof}

%\bibliographystyle{abbrv}
%\bibliography{/home/delmas/cermics/Bibliographie/delmas}
%\bibliography{biblio-feller}

\begin{thebibliography}{10}

\bibitem{ad:falp}
R.~ABRAHAM and J.-F. DELMAS.
\newblock Fragmentation associated to {L}évy processes.
\newblock {\em Preprint CERMICS}, 2005.

\bibitem{ap:sac}
D.~ALDOUS and J.~PITMAN.
\newblock The standard additive coalescent.
\newblock {\em Ann. Probab.}, 26(4):1703--1726, 1998.

\bibitem{b:pl}
J.~BERTOIN.
\newblock {\em L\'evy processes}.
\newblock Cambridge University Press, Cambridge, 1996.

\bibitem{b:fpcbm}
J.~BERTOIN.
\newblock A fragmentation process connected to {B}rownian motion.
\newblock {\em Probab. Th. Rel. Fields}, 117:289--301, 2000.

\bibitem{b:ssf}
J.~BERTOIN.
\newblock Self-similar fragmentations.
\newblock {\em Ann. Inst. Henri Poincar\'e}, 38(3):319--340, 2000.

\bibitem{b:smssf}
J.~BERTOIN.
\newblock On small masses in self-similar fragmentations.
\newblock {\em Stoch. Process. and Appl.}, 109(1):13--22, 2004.

\bibitem{b:rfcp}
J.~BERTOIN.
\newblock {\em Random fragmentation and coagulation processes}.
\newblock To appear, 2006.

\bibitem{d:fhalp}
J.-F. DELMAS.
\newblock Fragmentation at height associated to {L}évy processes.
\newblock {\em Preprint CERMICS}, 2006.

\bibitem{dlg:rtlpsbp}
T.~DUQUESNE and J.-F. LE~GALL.
\newblock {\em Random trees, {L}évy processes and spatial branching processes},
  volume 281.
\newblock Astérisque, 2002.

\bibitem{h:rfdsf}
B.~HAAS.
\newblock Regularity of formation of dust in self-similar fragmentations.
\newblock {\em Ann. Inst. Henri Poincar\'e}, 40(4):411--438, 2004.

\bibitem{lglj:bplplfss}
J.-F. LE~GALL and Y.~LE~JAN.
\newblock Branching processes in {L}évy processes: Laplace functionals of snake
  and superprocesses.
\newblock {\em Ann. Probab.}, 26:1407--1432, 1998.

\bibitem{lglj:bplpep}
J.-F. LE~GALL and Y.~LE~JAN.
\newblock Branching processes in {L}évy processes: The exploration process.
\newblock {\em Ann. Probab.}, 26:213--252, 1998.

\bibitem{m:sfdfstsn}
G.~MIERMONT.
\newblock Self-similar fragmentations derived from the stable tree {II}:
  splitting at nodes.
\newblock {\em Probab. Th. Rel. Fields}, 131(3):341--375, 2005.

\end{thebibliography}
\newcommand{\sortnoop}[1]{}

\end{document}